\documentstyle[twoside]{article} 
\newcommand{\es}{\emptyset}
\newcommand{\ci}{\subseteq}
\renewcommand{\u}{\cup}
\renewcommand{\i}{\cap}
\newcommand{\Lgra}{\Longrightarrow}
\renewcommand{\a}{\alpha}
\renewcommand{\b}{\beta}
\renewcommand{\d}{\delta}
\newcommand{\e}{\varepsilon}
\newcommand{\s}{\sigma}
\newcommand{\pf}{\noindent {\bf Proof :} }
\newcommand{\ov}{\overline}
\newcommand{\iy}{\infty}
\newcommand{\qed}{\hfill\rule{3mm}{3mm}}

\newcounter{cnt1}
\newcounter{cnt2}
\newcommand{\blr}{\begin{list}{$(\roman{cnt1})$}
	{\usecounter{cnt1} \setlength{\topsep}{0pt}
		\setlength{\itemsep}{0pt}}}
\newcommand{\bla}{\begin{list}{$($\alph{cnt2}$)$}
	{\usecounter{cnt2} \setlength{\topsep}{0pt}
		\setlength{\itemsep}{0pt}}}
\newcommand{\el}{\end{list}}

\newtheorem{thm}{Theorem}
\newtheorem{lem}[thm]{Lemma}
\newtheorem{cor}[thm]{Corollary}
\newtheorem{ex}{Example}
\newtheorem{Q}{Question}
\newtheorem{Def}{Definition}
\newtheorem{prop}[thm]{Proposition}
\newtheorem{rem}{Remark}
\newcommand{\Rem}{\begin{rem} \rm}
\newcommand{\bdfn}{\begin{Def} \rm}
\newcommand{\edfn}{\end{Def}}

\title{\bf The Mazur Intersection Property and Farthest
Points}
\author{Pradipta Bandyopadhyay\thanks{Stat--Math Division,
Indian Statistical Institute, 203, B.\ T.\ Road, Calcutta 700
035, {\sc India}, e-mail~: pradipta@isical.ernet.in
\newline 
{\bf AMS(MOS) Subject Classification (1980) :} 46B20, 46B22
}}
\date{}
\begin{document}

\maketitle
\begin{abstract} 
K.\ S.\ Lau had shown that a reflexive Banach space has the
Mazur Intersection Property (MIP) if and only if every closed
bounded convex set is the closed convex hull of its farthest
points.

In this work, we show that in general this latter property is
equivalent to a property stronger than the MIP. As corollaries,
we recapture the result of Lau and characterize the w*-MIP in
dual of RNP spaces.
\end{abstract}

\section*{Introduction}
We work with only real Banach spaces. The notations are
standard. Any unexplained terminology can be found in
\cite{Br,DU}.

\bdfn 
For a closed and bounded set $K$ in a Banach space $X$, define
\blr
\item $r_K(x) = \sup \{\|z-x\| : z\in K\}$, $x\in X$. $r_K$,
called the {\em farthest distance map}, is a Lipschitz
continuous convex function.
\item $Q_K(x) = \{z \in K : \|z-x\| = r_K(x)\}$, $x\in X$, the
set of points {\em farthest }\/ from $x$.
\item $D(K) = \{x \in X : Q_K(x) \neq \es\}$.
\item $b(K) = \u\{Q_K(x) : x\in D(K)\}$ is the set of {\em
farthest points}\/ of $K$.
\item for $x \in X$ and $\a>0$, let $C(K, x, \a) = \{z \in K :
\|z-x\| > r_K(x) - \a\}$. $C(K, x, \a)$ will be called {\em a
crescent of $K$ determined by $x$ and $\a$}.
\el

Call a closed and bounded set $K$ {\em densely remotal}\/ if
$D(K)$ is norm dense in $X$ and {\em almost remotal}\/ if $D(K)$
is generic, i.e., contains a dense $G_{\d}$ in $X$. $K$ {\em has
the Property $(R)$}\/ if any crescent of $K$ contains a farthest
point of $K$.
\edfn

\bdfn 
Call a set $K\ci X$ {\em admissible}\/ if it is the intersection
of closed balls containing it. Let $F$ be a norming subspace of
$X^*$ (i.e., $\|x\| = \sup \hat{x} (B(F))$, for all $x\in X$).
Let us denote the $\s(X,F)$ topology on $X$ simply by $\s$. Then
$B(X)$, the closed unit ball of $X$, is $\s$-closed and hence
any admissible set is $\s$-closed, bounded (in norm) and convex.
We say that {\em $X$ has $F$-MIP if}\/ the converse holds, i.e.,
{\em every bounded, $\s$-closed convex set in $X$ is
admissible}. This property was introduced in \cite{B2}
generalizing the Mazur Intersection Property (MIP) (i.e., when
$F=X^*$) and the w*-MIP (i.e., when $X=Y^*$ and $F=Y$) (see
\cite{B2,GGS} for various characterizations and related results).
\edfn

K.\ S.\ Lau \cite{L} had shown that in a reflexive space the MIP
is equivalent to the following~:
\begin{quote}
\em Every closed bounded convex set is the closed convex hull of
its farthest points.
\end{quote}

In this work, we show that in a Banach space $X$, every
$\s$-closed bounded convex set is the $\s$-closed convex hull of
its farthest points if and only if $X$ has the $F$-MIP and every
$\s$-closed bounded convex set in $X$ has the Property $(R)$. As
corollaries, we recapture Lau's result and characterize the
w*-MIP in dual of spaces with the Radon-Nikod\'{y}m Property
(RNP) using a result of Deville and Zizler \cite{DZ}.

Some of the results presented here was first observed in the
author's Ph.~D.\ Thesis \cite{PBTH} written under the
supervision of Prof.\ A.\ K.\ Roy.

\section{Main Results}
The following Lemma is an easy consequence of the fact that the
class of admissible sets is closed under arbitrary intersection.
\begin{lem} 
\label{IntMIP} Let $X$ be a Banach space and $F$ be a norming
subspace of $X^*$. Let $K \ci X$ be $\s$-closed, bounded and
convex. If for all $\l>0$, the set $K_{\l}=\ov{K+\l B(X)}^{\s}$
is admissible, then so is $K$.
\end{lem}
\Rem 
\label{r1} It follows that given any $\s$-closed, bounded convex
inadmissible set $K \ci X$, there exists an inadmissible set of
the form $K_{\l}$, which has non-empty norm interior. This
improves \cite[Lemma 3.2]{L}.
\end{rem} 
\begin{Q} 
Is the converse true~?
\end{Q}
\begin{thm} 
\label{far} Let $X$ be a Banach space and $F$ be a norming
subspace of $X^*$. The following statements are equivalent~:
\bla
\item Any $\s$-closed, bounded convex set in $X$ is the
$\s$-closed convex hull of its farthest points.
\item $(i)$ $X$ has the $F$-MIP, and 

$(ii)$ Any $\s$-closed, bounded convex set in $X$ has the
Property $(R)$.
\el 
\end{thm}

\pf $(a) \Lgra (b)$ We first prove $(i)$ following Lau's lead
\cite{L}. Suppose $X$ lacks the $F$-MIP. By Remark~\ref{r1},
there exists a $\s$-closed, bounded convex set $K \ci X$ with
int$(K) \neq \es$ that is not admissible.

Let $M = \i \{ B : B$ closed ball containing $K \}$ and let
$x_o \in M \setminus K$ and $y_o \in \mbox{int}(K) \ci \mbox{int}
(M)$. Since $M \setminus K$ is open in $M$, there exists $0 < \l
<1$ such that $z_o = \l x_o + (1-\l)y_o \in M \setminus K$. Note
that $z_o \in \mbox{int}(M)$, and hence, so is any point of the
form
\begin{quote}
$(*)$ \hfill $\a z_o + (1-\a)x, \quad \a \in (0, 1], \quad x \in
K$. \hfill~
\end{quote} 

Let $K_1 = \mbox{co}(K \u \{z_o\})$. Then $K_1$ is $\s$-closed,
bounded and convex. The proof will be complete once we show that
$b(K_1) \ci K$.

Let $x \in X$. Then $B = \{z \in X : \|z-x\| \leq r_K (x) \}$ is
a closed ball containing $K$, and so contains $M$. Since each
point of the form $(*)$ is in $\mbox{int} (M)$, it is in int$(B)$,
i.e., its distance from $x$ is strictly less than $r_K(x)$.
Note that $r_K(x) \leq r_{K_1}(x) \leq r_M(x) = r_K(x)$. Thus,
$Q_{K_1}(x) \ci K$. Since $x \in X$ was arbitrary, $b(K_1) \ci
K$, contradicting (a).

Now, if there exists a $\s$-closed, bounded convex set $K$ and a
crescent of $K$ that is disjoint from $b(K)$, then it is
disjoint from $\ov{\rm co}^{\s} b(K)$ as well. Hence $\ov{\rm
co}^{\s} b(K) \neq K$. This proves $(ii)$.

$(b) \Lgra (a)$ Let $K$ be a $\s$-closed, bounded convex set in
$X$. Let $L = \ov{\rm co}^{\s} (b(K))$. Clearly, $L \ci K$.
Suppose there exists $x \in K \setminus L$. Since $X$ has the
$F$-MIP, there exists a crescent $C$ of $K$ disjoint from $L$.
Since $K$ has the Property $(R)$, $C \i b(K) \neq \es$. But, of
course, $b(K) \ci L$.
\qed
\begin{lem} 
\label{open} Let $K \ci X$ be a bounded set. Let $x \in X$ and
$\a > 0$ be given. Then there exists $\e > 0$ such that for any
$y \in X$ with $\|x-y\| < \e$, there exists $\b > 0$ such that
$C(K, y, \b) \ci C(K, x, \a)$.
\end{lem}

\pf Take $0 < \e < \a/2$ and $0 < \b < \a - 2\e$.
\qed
\begin{prop} 
\label{R} Any densely remotal set has Property $(R)$.
\end{prop}

\pf If $D(K)$ is dense in $X$ and $C(K, x, \a)$ is any crescent
of $K$, then, by Lemma~\ref{open}, there exists $y \in D(K)$ and
$\b > 0$ such that $C(K, y, \b) \ci C(K, x, \a)$. Clearly,
$Q_K(y) \ci C(K, y, \b)$. Thus, $K$ has the Property $(R)$.
\qed

In the following Lemma, we collect some known results that
identify some important classes of sets with the Property $(R)$.
\begin{lem} 
\label{remotal} $(a)$ {\rm \cite[Theorem 2.3]{L}} Any weakly
compact set is almost remotal with respect to any equivalent
norm.

$(b)$ {\rm \cite[Proposition 3]{DZ}} If $X$ has the RNP, every
w*-compact set in $X^*$ is almost remotal with respect to any
equivalent dual norm.
\end{lem}
\Rem 
Proposition 1 of \cite{DZ} gives an example of a w*-compact
convex subset of $\ell^1$ that lacks farthest points. Thus, some
additional assumptions are necessary even for Property $(R)$.
\end{rem}

Combining this with Theorem~\ref{far}, we get
\begin{thm} 
\label{far2} If $X$ has the RNP, then $X^*$ has the w*-MIP if
and only if every w*-compact convex set in $X^*$ is the
w*-closed convex hull of its farthest points.
\end{thm}
\begin{cor} 
\label{far1} {\rm \cite{L}} If $X$ is reflexive, then $X$ has
the MIP if and only if every closed bounded convex set in $X$ is
the closed convex hull of its farthest points.
\end{cor}

In both the above cases, the additional assumption on $X$
already implies condition $(b) (ii)$ of Theorem~\ref{far}, and
hence the equivalence. However, in the corollary below, this is
not so direct.

Recall (from \cite{GK}) that a space has the IP$_{f,\,\iy}$ if
every family of closed balls with empty intersection contain a
finite subfamily with empty intersection. For example, any space
$X$ that is 1-complemented in a dual space, in particular any
dual space, has the IP$_{f,\,\iy}$.
\begin{cor} 
\label{IP} If $X$ has the IP$_{f,\,\iy}$, then $X$ has the MIP
if and only if every closed bounded convex set in $X$ is the
closed convex hull of its farthest points.
\end{cor}

\pf Sufficiency follows from Theorem~\ref{far}. Conversely,
if $X$ has both the MIP and IP$_{f,\,\iy}$ then it must be
reflexive \cite[Theorem~VIII.5]{GK}.
\qed
\Rem 
(a) From the known characterizations it is easily seen that if
$X$ has the MIP then $X^{**}$ has the w*-MIP. And it is a long
standing open question whether then $X$ is also Asplund, or
equivalently, $X^*$ has the RNP (see \cite{GGS}). If the answer
to this question is yes, then by Theorem~\ref{far2}, so is that
to the following
\begin{Q} 
If $X$ has the MIP, is every w*-compact convex set in $X^{**}$
the w*-closed convex hull of its farthest points~?
\end{Q}

(b) The proof of $(b) \Lgra (a)$ in Theorem~\ref{far}, combined
with Lemma~\ref{remotal}, also shows that if every weakly
compact or compact convex set is admissible then each such set
is the closed convex hull of its farthest points. Is the
converse true~? Clearly, a similar proof will not work unless
the space is reflexive or finite dimensional. Now, can the
specific form of $K_{\l}$, as in Lemma~\ref{IntMIP}, be utilized
to prove the converse~?

(c) We see from the last three results that in some cases the
condition (a) of Theorem~\ref{far} becomes equivalent to the
$F$-MIP. Is this generally true~? We do not know the answer for
arbitrary $F$. For example, we do not know whether
Theorem~\ref{far2} holds even without the RNP. However, the
answer is negative for $F=X^*$ as the following example shows.
\end{rem}
\begin{ex} 
There is a space $X$ with a Fr\'{e}chet differentiable norm, and
hence with MIP, and a closed bounded convex set $K$ in $X$ that
lacks farthest points.
\end{ex}

\pf Notice that if the norm is strictly convex (respectively,
locally uniformly convex), any farthest point of a closed
bounded convex set is also an extreme (resp.\ denting) point.
So, if every closed bounded convex set admits farthest points,
then the space must necessarily have the Krein-Milman Property
(KMP) (resp.\ the RNP) (see \cite{Br} or \cite{DU} for details).
However, the space $c_o$, which does not have the KMP, admits a
strictly convex Fr\'{e}chet differentiable norm (see e.g.,
\cite{Di}).
\qed
\Rem 
Observe that since $c_o$ is Asplund, Theorem~\ref{far2} shows
that when equipped with the bidual norm of the above, every
w*-compact convex set in $\ell^{\iy}$ is the w*-closed convex
hull of its farthest points. Thus, there is closed bounded
convex set $K \ci c_o$, such that no farthest point of
$\widetilde{K}$, the w*-closure of the canonical image of $K$ in
$X^{**}$, comes from $K$.
\end{rem}
 
This gives rise to two very natural questions.
\begin{Q} 
$(a)$ If $X$ has both RNP and MIP, does every closed bounded
convex set in $X$ have the Property $(R)$~?

$(b)$ Can the condition $(b) (ii)$ of Theorem~\ref{far} be
replaced by the weaker condition that every $\s$-closed, bounded
convex set in $X$ has farthest points~?
\end{Q}
\Rem 
The example in Proposition 1 of \cite{DZ} shows also that RNP
alone is not enough to ensure even existence of farthest points.
But then, the space there does not even have the w*-MIP. In
fact, the set under consideration itself is not admissible.
Indeed, the answer to the first question is likely to be
affirmative. As evidence, observe that
\bla
\item if $z \in K$ is farthest from $x \in X$, then it is
nearest from any point on the line $[x, \, z]$ extended beyond
$z$ (this is just triangle inequality, and was observed in
\cite{SF}), and
\item if $X$ has both RNP and MIP, any crescent of any closed
bounded convex set $K$ in $X$ contains a {\em nearest point}\/
of $K$. This is because in spaces with RNP, any closed bounded
convex set is the closed convex hull of its nearest points
\cite[Theorem 8.3]{BF}.
\el
One thus possibly needs to characterize farthest points among
nearest points.
\end{rem}

\end{document}